\documentclass[numreferences]{kluwer}
\usepackage{amssymb}
\usepackage{graphics}

\newcommand{\ds}{\displaystyle}
\newcommand{\bc}{\begin{equation}\begin{array}{l}}
\newcommand{\ec}{\end{array}\end{equation}}
\newcommand{\ason}{\renewcommand{\arraystretch}{1.5}}
\newcommand{\asoff}{\renewcommand{\arraystretch}{1.0}}
\newcommand{\nn}{\nonumber}

\begin{document}
\begin{article}
\begin{opening}

\title{High Order Multistep Methods with Improved Phase-Lag Characteristics for the Integration of the\\ Schr\"{o}dinger Equation}

\author{D.S. \surname{Vlachos}\thanks{e-mail: dvlachos@uop.gr}}

\author{Z.A. \surname{Anastassi}\thanks{e-mail: zackanas@uop.gr}}

\author{T.E. \surname{Simos}\thanks{Highly Cited Researcher, Active Member of the European Academy of Sciences and
Arts. Corresponding Member of the European Academy of Sciences Corresponding Member of European Academy of Arts, Sciences and Humanities, Please use the following address for all correspondence: Dr. T.E. Simos, 10 Konitsis Street, Amfithea - Paleon Faliron, GR-175 64 Athens, Greece, Tel: 0030 210 94 20 091, e-mail: tsimos.conf@gmail.com, tsimos@mail.ariadne-t.gr}}

\institute{Laboratory of Computational Sciences, Department of Computer Science and Technology, Faculty of Sciences
and Technology, University of Peloponnese, GR-221 00 Tripolis,
Greece}

\runningtitle{High Order Multistep Methods with Improved Phase-Lag Characteristics}
\runningauthor{D.S. Vlachos, Z.A. Anastassi, T.E. Simos}

\begin{abstract}
In this work we introduce a new family of twelve-step linear multistep methods for the integration of the Schr\"odinger equation. The new methods are constructed by adopting a new methodology which improves the phase lag characteristics by vanishing both the phase lag function and its first derivatives at a specific frequency. This results in decreasing the sensitivity of the integration method on the estimated frequency of the problem. The efficiency of the new family of methods is proved via error analysis and numerical applications.
\end{abstract}

\keywords{Numerical solution, Schr\"odinger equation, multistep methods, hybrid methods, P-stability, phase-lag, phase-fitted}

\classification{PACS}{02.60, 02.70.Bf, 95.10.Ce, 95.10.Eg, 95.75.Pq}
\end{opening}

\section{Introduction}
\label{intro}
The numerical integration of systems of ordinary differential equations with
oscillatory solutions has been a subject of research during the past decades.
This type of ODEs is often met in real problems, like the Schr\"{o}dinger equation
and the N-body problem. For problems having highly oscillatory solutions, standard
methods with unspecialized use can require a huge number of steps to track the oscillations. One way to obtain a more efficient integration process is to construct numerical
methods with an increased algebraic order, although the implementation of high algebraic order meets several difficulties \cite{quinlan_arxiv_astro_ph_9901136}.

On the other hand, there are some special techniques for optimizing numerical methods. Trigonometrical fitting and phase-fitting are some of them, producing methods with
variable coefficients, which depend on $v = \omega h$, where $\omega$ is the dominant frequency of the problem and $h$ is the step length of integration. More precisely, the coefficients of a general linear method are found from
the requirement that it integrates exactly powers up to degree $p+1$. For problems having oscillatory
solutions, more efficient methods are obtained when they are exact for every linear
combination of functions from the reference set
\begin{equation}
\{1, x, \ldots , x^K , e^{\pm \mu x},\ldots , x^P e^{\pm \mu x}\}\label{equ_exp_fit}
\end{equation}
This technique is known as exponential (or trigonometric if $\mu=i\omega$) fitting and has a long history \cite{gautschi_NM_3_381_61}, \cite{lyche_NM_19_65_72}. The
set (\ref{equ_exp_fit}) is characterized by two integer parameters, $K$ and $P$ . The set in which there
is no classical component is identified by $K =-1$ while the set in which there is no
exponential fitting component (the classical case) is identified by $P =-1$. Parameter
$P$ will be called the level of tuning. An important property of exponential fitted algorithms is that they tend to the corresponding classical ones when the involved frequencies tend to zero, a fact which allows us to say that exponential fitting represents a natural extension of the classical polynomial fitting. The examination of the convergence
of exponential fitted multistep methods is included in Lyche's theory \cite{lyche_NM_19_65_72}. There is a large number of significant methods presented with high practical importance that have been presented in the bibliography. The general theory is presented in detail in \cite{ixaru_Book_EF_KAP_04}.

Considering the accuracy of a method, when solving oscillatory problems,
it is more appropriate to work with the phase-lag, rather than the principal local truncation error. We mention
the pioneering paper of Brusa and Nigro \cite{brusa_IJNME_15_685_80}, in which the phase-lag property was
introduced. This is actually another type of a truncation error, i.e. the angle between
the analytical solution and the numerical solution. On the other hand, exponential fitting is accurate only when a good estimate of the dominant frequency of the solution
is known in advance. This means that in practice, if a small change in the dominant frequency is introduced, the efficiency of the method can be dramatically altered. It is well known that for equations similar to the harmonic oscillator the most efficient exponential fitted methods are
those with the highest tuning level. A lot of significant work has been made during the last years in this field, mainly focusing for obvious reasons in the solution of the Schr\"{o}dinger equation (see for example \cite{ix78}-\cite{jnaiam3_11}).

In this paper we present a new family of methods based on the 12-step linear multistep method of Quinlan and Tremaine \cite{quinlan_AJ_100_1694_90}. The new methods are constructed by vanishing the phase-lag function and its first derivatives at a predefined frequency. Error analysis and numerical experiments show that the new methods exhibit improved characteristics concerning the solution of the time-independent Schr\"odinger equation. The paper is organized as follows: In section 2 the general theory of the new methodology is presented. In section 3 the new methods are described in detail. In section 4 the stability properties of the new methods are investigated. Section 5 presents the results from the numerical experiments and finally, conclusions are drawn in section 6.

\section{Phase-lag analysis of symmetric multistep methods}
Consider the differential equations
\begin{equation}
\frac{d^2y(t)}{dt^2}=f(t,y),\;y(t_0)=y_0,\;y'(t_0)=y'_0
\label{equ_ref_diff}
\end{equation}
and the linear multistep methods
\begin{equation}
\sum_{j=0}^{J}a_jy_{n+j}=h^2\sum_{j=0}^{J}b_jf_{n+j}
\label{equ_ref_meth}
\end{equation}
where $y_{n+j}=y(t_0+(n+j)h)$, $f_{n+j}=f(t_0+(n+j)h,y(t_0+(n+j)h))$ and $h$ is the step size of the method. With the method (\ref{equ_ref_meth}), we associate the following functional
\begin{equation}
L(h,a,b,y(t))=\sum_{j=0}^Ja_jy(t+j\cdot h)-h^2\sum_{j=0}^Jb_jy''(t+j\cdot h)
\end{equation}
where $a,b$ are the vectors of coefficients $a_j$ and $b_j$ respectively, and $y(t)$ is an arbitrary function. The algebraic order of the method (\ref{equ_ref_meth}) is $p$, if
\begin{equation}
L(h,a,b,y(t))=C_{p+2}h^{p+2}y^{(p+2)}(t)+O(h^{p+3})
\label{equ_ref_error}
\end{equation}
The coefficients $C_q$ are given
\begin{eqnarray}
C_0=&&\sum_{j=0}^J a_j \nonumber \\
C_1=&&\sum_{j=0}^J j\cdot a_j \nonumber \\
c_q=&&\frac{1}{q!}\sum_{j=0}^Jj^q\cdot a_j -\frac{1}{(q-2)!}\sum_{j=0}^Jj^{q-2}b_j
\end{eqnarray}
The principal local truncation error (PLTE) is the leading term of (\ref{equ_ref_error})
\begin{equation}
 plte=C_{p+2}h^{p+2}y^{(p+2)}(t)
\end{equation}
The following assumptions will be considered in the rest of the paper:
\begin{enumerate}
\item $a_J=1$, since we can always divide the coefficients of (\ref{equ_ref_meth}) with $a_J$.
\item $|a_0|+|b_0|\neq 0$, since otherwise we can assume that $J=J-1$.
\item $\sum_{j=0}^J |b_j| \neq 0$, since otherwise the solution of (\ref{equ_ref_meth}) would be independent of (\ref{equ_ref_diff}).
\item The method (\ref{equ_ref_meth}) is at least of order one.
\item The method (\ref{equ_ref_meth}) is zero stable, which means that the roots of the polynomial
\begin{equation}
p(z)=\sum_{j=0}^Ja_jz^j
\end{equation}
all lie in the unit disc, and those that lie on the unit circle have multiplicity one.
\item The method (\ref{equ_ref_meth}) is symmetric, which means that
\begin{equation}
a_j=a_{J-j},\;b_j=b_{J-j},\;j=0(1)J
\end{equation}
It is easily proved that both the order of the method and the step number $J$ are even numbers \cite{lambert_JIMA_18_189_76}.
\end{enumerate}
Consider now the test problems
\begin{equation}
y''(t)=-\omega ^2 y(t)
\label{equ_ref_equ}
\end{equation}
where $\omega$ is a constant. The numerical solution of (\ref{equ_ref_equ}) by applying method (\ref{equ_ref_meth}) is described by the difference equation
\begin{equation}
\sum_{j=1}^{J/2} A_j(s^2)(y_{n+j}+y_{n-j})+A_0(s^2)y_n=0
\end{equation}
with
\begin{equation}
A_j(s^2)=a_{\frac{J}{2}-j}+s^2\cdot b_{\frac{J}{2}-j}
\end{equation}
and $s=\omega h$. The characteristic equation is then given by
\begin{equation}
\sum_{j=1}^{J/2} A_j(s^2)(z^j+z^{-j})+A_0(s^2)=0
\label{equ_char_equ}
\end{equation}
and the interval of periodicity $(0,s_0^2)$ is then defined such that for $s\in (0,s_0)$ the roots of (\ref{equ_char_equ}) are of the form
\begin{equation}
z_1=e^{i\lambda (s)},\;z_2=e^{-i\lambda (s)},\;|z_j|\leq 1,\;3\leq j\leq J
\end{equation}
where $\lambda (s)$ is a real function of $s$. The phase-lag $PL$ of the method (\ref{equ_ref_meth}) is then defined
\begin{equation}
PL=s-\lambda (s)
\end{equation}
and is of order $q$ if
\begin{equation}
PL=c\cdot s^{q+1}+O(s^{q+3})
\end{equation}
In general, the coefficients of the method (\ref{equ_ref_meth}) depend on some parameter $v$, thus the coefficients $A_j$ are functions of both $s^2$ and $v$. The following theorem was proved by Simos and Williams \cite{simos_CC_23_513_99}:
%\begin{theorem}
For the symmetric method (\ref{equ_ref_equ}) the phase-lag is given
\begin{equation}
PL(s,v)=\frac{2\sum_{j=1}^{J/2}A_j(s^2,v)\cdot cos(j\cdot s) +A_0(s^2,v)}{2\sum_{j=1}^{J/2}j^2A_j(s^2,v)}
\end{equation}
%\end{theorem}
We are now in position to describe the new methodology. In order to efficiently integrate oscillatory problems, it is a good practice to calculate the coefficients of the numerical method by forcing the phase lag to be zero at a specific frequency. But, since the appropriate frequency is problem dependent and in general is not always known, we may assume that we have an error in the frequency estimation. It would be of great importance to force the phase-lag to be insensitive to this error. Thus, beyond the vanishing of the phase-lag, we also force its first derivatives to be zero.

\section{Construction of the new methods}
\subsection{Classical Method}
The family of new methods is based on the 12-step linear multistep method of Quinlan and Tremaine \cite{quinlan_AJ_100_1694_90} which is of the form (\ref{equ_ref_meth}) with coefficients

\ason
\begin{equation}
\begin{array}{cccc}
a_0=1 & a_1=-2 & a_2=2 & a_3=-1\\
a_4=0 & a_5=0 & a_6=0\\
b_0=0 & b_1=\frac{90987349}{53222400} & b_2=-\frac{-114798419}{26611200} & b_3=\frac{270875723}{17740800}\\
b_4=-\frac{-67855831}{2217600} & b_5=\frac{50277247}{985600} & b_6=\frac{-253491379}{4435200}
\end{array}
\label{equ_base_meth}
\end{equation}
\asoff

The PLTE of the method is given:
\begin{equation}
 plte=\left\{\frac{16301796103 y^{(14)} \
h^{14}}{290594304000}+O\left(h^{16}\right)\right\}
\end{equation}
\subsection{New Methods using Phase Fitting}
The methods that are constructed are named as \textit{PF-Di}, where:
\begin{itemize}
 \item \textit{PF-D0}: the phase lag function is zero at the frequency $v=\omega *h$.
\item \textit{PF-D1}: the phase lag function and its first derivative are zero at the frequency $v=\omega *h$.
\item \textit{PF-D2}: the phase lag function and its first and second derivatives are zero at the frequency $v=\omega *h$.
\item \textit{PF-D3}: the phase lag function and its first, second and third derivatives are zero at the frequency $v=\omega *h$.
\item \textit{PF-D4}: the phase lag function and its first, second, third and fourth derivatives are zero at the frequency $v=\omega *h$.
\item \textit{PF-D5}: the phase lag function and its first, second, third, fourth and fifth derivatives are zero at the frequency $v=\omega *h$.
\end{itemize}
The coefficients of the methods in the form:
\begin{eqnarray}
b^i_1=&&\frac{b_{1,num}^i}{b_{1,denum}^i} \nonumber \\
b^i_2=&&\frac{b_{2,num}^i}{b_{2,denum}^i} \nonumber \\
b^i_3=&&\frac{b_{3,num}^i}{b_{3,denum}^i} \nonumber \\
b^i_4=&&\frac{b_{4,num}^i}{b_{4,denum}^i} \nonumber \\
b^i_5=&&\frac{b_{5,num}^i}{b_{5,denum}^i} \nonumber \\
b^i_6=&&\frac{b_{6,num}^i}{b_{6,denum}^i} \nonumber
\end{eqnarray}
where the coefficients $b^i$ correspond to the method \textit{PF-Di}. Since for small values of $v$, the above formulae are subject to heavy cancelations, the Taylor expansions of the coefficients have been calculated as $b_T^i$. The exact formula of all coefficients are given in appendix.

The principal local truncation errors of the methods are given by

$$\ds plte_i = \frac{16301796103}{290594304000} coeff_i,$$ where

\begin{equation}
\begin{array}{ll}
coeff_0 = &{\left(y^{(12)} \omega ^2+y^{(14)}\right) \
h^{14}}
\nonumber \\
coeff_1 = &{ \left(y^{(10)} \omega ^4+2 y^{(12)} \omega \
^2+y^{(14)}\right) h^{14}}
\nonumber \\
coeff_2 = &{ \left(y^{(8)} \omega ^6+3 y^{(10)} \omega \
^4+3 y^{(12)} \omega ^2+y^{(14)}\right) \
h^{14}}
\nonumber \\
coeff_3 = &{ \left(y^{(6)} \omega ^8+4 y^{(8)} \omega ^6+6 \
y^{(10)} \omega ^4+4 y^{(12)} \omega ^2+y^{(14)}\right) h^{14}}
\nonumber \\
coeff_4 = &{ \left(y^{(4)} \omega ^{10}+5 y^{(6)} \omega \
^8+10 y^{(8)} \omega ^6+10 y^{(10)} \omega ^4\right.}\\
&\nn {\left.+5 y^{(12)} \omega ^2+y^{(14)}\right) h^{14}}
\nonumber \\
coeff_5 = &{ \left(y^{(2)} \omega ^{12}+6 y^{(4)} \omega \
^{10}+15 y^{(6)} \omega ^8+20 y^{(8)} \omega ^6\right.}\\
&\nn {\left.+15 y^{(10)} \omega ^4+6 y^{(12)} \
\omega ^2+y^{(14)}\right) h^{14}}
\ec

\section{Stability Analysis}
The stability of the new methods is studied by considering the test equation
\begin{equation}
 \frac{d^2y(t)}{dt^2}=-\sigma ^2 y(t)
\end{equation}
and the linear multistep method (\ref{equ_ref_meth}) for the numerical solution.
In the above equation $\sigma \neq \omega$ ($\omega$ is the frequency at which the phase-lag function and its derivatives vanish). By setting $s=\sigma h$ and $v=\omega h$, we get for the characteristic equation of the applied method
\begin{equation}
 \sum_{j=1}^{J/2} A_j(s^2,v)(z^j+z^{-j})+A_0(s^2,v)=0
\end{equation}
where
\begin{equation}
 A_j(s^2,v)=a_{\frac{J}{2}-j}(v)+s^2\cdot b_{\frac{J}{2}-j}(v)
\end{equation}
The motivation of the above analysis is straightforward: Although the coefficients of the method (\ref{equ_ref_meth}) are designed in a way that the phase-lag and its first derivatives vanish in the frequency $\omega$, the frequency $\omega$ itself is unknown and only an estimation can be made. Thus, if the correct frequency of the problem is $\sigma$ we have to check if the method is stable, that is if the roots of the characteristic equation lie in the unit disk. For this reason we draw in the $s-v$ plane the areas in which the method is stable. Figure \ref{fig:1} shows the stability region for the classical method and figure \ref{fig:2} for the six methods (the phase fitted one and those with first,second, third, fourth and fifth phase lag derivative elimination). Note here that the $s$-axis corresponds to the real frequency while the $v$-axis corresponds to the estimated frequency used to construct the parameters of the method.

\section{Numerical Results}
The radial Schr\"odinger equation is given by:
\begin{equation}
 y''(x)=\left(\frac{l(l+1)}{x^2}+V(x)-E)y(x)\right)
\label{def_schr_equ}
\end{equation}
where $\frac{l(l+1)}{x^2}$ is the centrifugal potential, $V(x)$ is the potential, $E$ is the Energy and $W(x)=\frac{l(l+1)}{x^2}+V(x)$ is the effective potential. It is valid that $lim_{x\rightarrow \infty}V(x)=0$ and therefore $lim_{x\rightarrow \infty}W(x)=0$.
We consider that $E>0$ and we divide the interval $[0,+\infty)$ into subintervals $[a_i,b_i]$ so that $W(x)$ can be considered constant inside each subinterval with value $\hat{W}_i$. The problem (\ref{def_schr_equ}) can be expressed now by the equations
\begin{equation}
 y''_i=(\hat{W}_i-E)y_i
\end{equation}
whose solution are
\begin{equation}
 y_i(x)=\left(A_ie^{\sqrt{\hat{W}_i-E}x}+B_ie^{-\sqrt{\hat{W}_i-E}x}\right)
\end{equation}
with $A_i,B_i \in R$.
We will integrate problem (\ref{def_schr_equ}) with $l = 0$ at the interval $[0, 15]$ using the well
known Woods-Saxon potential:
\begin{equation}
 V(x)=\frac{u_0}{1+q}+\frac{U_1q}{(1+q)^2}, \quad q=e^{\frac{x-x_0}{a}}
\end{equation}
where $u_0=-50$, $a=0.6$, $x_0=7$, $u_1=-\frac{u_0}{a}$ and with boundary condition $y(0)=0$. The potential $V(x)$ decays more quickly than $\frac{l(l+1)}{x^2}$ , so for large x (asymptotic region) the Schr\"odinger equation (\ref{def_schr_equ}) becomes
\begin{equation}
 y''(x)=\left(\frac{l(l+1)}{x^2}-E)y(x)\right)
\end{equation}
The last equation has two linearly independent solutions $kxj_l(kx)$ and
$kxn_l(kx)$, where $j_l$ and $n_l$ are the spherical \textit{Bessel} and \textit{Neumann} functions and $k=\sqrt{\frac{l(l+1)}{x^2}-E}$.
When $x\to \infty$ the solution takes the asymptotic form
\begin{eqnarray}
 y(x) \sim && A k x j_l (k x) - B k x n_l (k x) \nonumber \\
\sim && D[sin(k x - \pi \frac{l}{2}) + tan(\delta_l ) \cos (k x - \pi \frac{l}{2})],
\end{eqnarray}
where $\delta _l$ is called the \textit{scattering phase shift} and it is given by the following expression:
\begin{equation}
 tan(\delta _l)=\frac{y(x_i)S(x_{i+1})-y(x_{i+1})S(x_i)}{y(x_{i+1}C(x_i)-y(x_i)C(x_{i+1})}
\end{equation}
where $S(x)=kxj_l(kx)$ and $C(x)=kxn_l(kx)$ and $x_i<x_{i+1}$                            and both belong to the asymptotic region. Given the energy, we approximate the phase shift, the
accurate value of which is $\frac{\pi}{2}$ for the above problem.
We will use three different values for the energy: i) $989.701916$, ii) $341.495874$
and iii) $163.215341$. As for the frequency $\omega$ we will use the suggestion of Ixaru
and Rizea \cite{ixaru_CPC_38_3329_85}:
\begin{equation}
 \omega=\left\{\begin{array}{ll}\sqrt{E-50}, & x\in [0,6.5] \\ \sqrt{E}, & x\in [6.5,15] \end{array} \right.
\end{equation}
The results are shown in figures \ref{fig:3}, \ref{fig:4} and \ref{fig:5}. It is clear that the accuracy increases as the number of the eliminated derivatives of the phase lag function increases.

\section{Conclusions}
We have presented a new family of 12-step symmetric multistep numerical methods with improved characteristics concerning the integration of the Schr\"odinger equation. The methods were constructed by adopting a new methodology which, except for the phase fitting at a predefined frequency, it eliminates the first derivatives of the phase lag function at the same frequency. The result is that the phase lag function becomes less sensitive on the frequency near the predefined one. This behavior compensates the fact that the exact frequency can only be estimated. Experimental results demonstrate this behavior by showing that the accuracy is increased as the number of the derivatives that are eliminated is increased.

\clearpage

\appendix
Method \textit{PF-D0}:
\bc
b_{1,num}^0=((-124184636 \cos (v)+70378348 \cos (2 v)\\
-24862148 \cos (3 v)+5153611 \cos (4 v)) v^2\\
+25 (3013169 v^2-16128 \cos (3
v)\\
+32256 \cos (4 v)-32256 \cos (5 v)+16128 \cos (6 v)))
\csc ^{10}(\frac{v}{2})
\nonumber \\
b_{2,num}^0=(v^2 (159588050 \cos (v)-85350160 \cos (2
v)\\
+16708985 \cos (3 v)-5153611 \cos (5 v))-16 (6496079
v^2\\
-252000 \cos (3 v)+504000 \cos (4 v)\\
-504000 \cos (5 v)+252000 \cos
(6 v))) \csc ^{10}(\frac{v}{2})
\nonumber \\
b_{3,num}^0=((-367257540 \cos (v)+183567900 \cos (2 v)\\
-16708985
\cos (4 v)+24862148 \cos (5 v)) v^2\\
+81 (3175117 v^2-224000 \cos
(3 v)+448000 \cos (4 v)\\
-448000 \cos (5 v)+224000 \cos (6 v))) \csc ^{10}(\frac{v}{2})
\nonumber \\
b_{4,num}^0=(30675810 \cos (v) v^2-42958788 v^2\\
+75
(161280-611893 v^2) \cos (3 v)\\
+140 (152411 v^2-172800
) \cos (4 v)\\
+(24192000-17594587 v^2) \cos (5
v)-12096000 \cos (6 v)) \csc ^{10}(\frac{v}{2})
\nonumber \\
b_{5,num}^0=(-61351620 \cos (2 v) v^2+85936557 v^2\\
+30
(6120959 v^2-1411200) \cos (3 v)\\
+25 (3386880-3191761
v^2) \cos (4 v)\\
+(62092318 v^2-84672000) \cos (5
v)+42336000 \cos (6 v)) \csc ^{10}(\frac{v}{2})
\nonumber \\
b_{6,num}^0=((-171873114 \cos (v)+171835152 \cos (2 v)\\
-257184477
\cos (3 v)+103937264 \cos (4 v)\\
-75329225 \cos (5 v)) v^2+50803200
(\cos (3 v)-2 \cos (4 v)\\
+2 \cos (5 v)-\cos (6 v))) \csc
^{10}(\frac{v}{2})
\nonumber \\
b_{1,denom}^0=206438400 v^2
\nonumber \\
b_{2,denom}^0=206438400 v^2
\nonumber \\
b_{3,denom}^0=206438400 v^2
\nonumber \\
b_{4,denom}^0=51609600 v^2
\nonumber \\
b_{5,denom}^0=103219200 v^2
\nonumber \\
b_{6,denom}^0=103219200 v^2
\nonumber
\ec

\bc
\ds b^0_{T,1}=\frac{90987349}{53222400}-\frac{16301796103 \
v^2}{290594304000}+\frac{2012122579 v^4}{581188608000}\\
\ds -\frac{48410309 \
v^6}{1184440320000}+\frac{991945331743 \
v^8}{2838385676206080000}-\ldots\\
\ds b^0_{T,2}=-\frac{114798419}{26611200}+\frac{16301796103 \
v^2}{29059430400}-\frac{2012122579 v^4}{58118860800}\\
\ds +\frac{48410309 \
v^6}{118444032000}-\frac{991945331743 \
v^8}{283838567620608000}+\ldots
\nonumber \\
\ds b^0_{T,3}=\frac{270875723}{17740800}-\frac{16301796103 \
v^2}{6457651200}+\frac{2012122579 v^4}{12915302400}\\
\ds -\frac{48410309 \
v^6}{26320896000}+\frac{991945331743 \
v^8}{63075237249024000}-\ldots
\nonumber \\
\ds b^0_{T,4}=-\frac{67855831}{2217600}+\frac{16301796103 \
v^2}{2421619200}-\frac{2012122579 v^4}{4843238400}\\
\ds +\frac{48410309 \
v^6}{9870336000}-\frac{991945331743 \
v^8}{23653213968384000}+\ldots
\nonumber \\
\ds b^0_{T,5}=\frac{50277247}{985600}-\frac{16301796103 \
v^2}{1383782400}+\frac{2012122579 v^4}{2767564800}\\
\ds -\frac{48410309 \
v^6}{5640192000}+\frac{991945331743 \
v^8}{13516122267648000}-\ldots
\nonumber \\
\ds b^0_{T,6}=-\frac{253491379}{4435200}+\frac{16301796103 \
v^2}{1153152000}-\frac{2012122579 v^4}{2306304000}\\
\ds +\frac{48410309 \
v^6}{4700160000}-\frac{991945331743 \
v^8}{11263435223040000}+\ldots\\
\nonumber
\ec

Method \textit{PF-D1}:
\bc
b_{1,num}^1=\csc ^{10}(\frac{v}{2}) \sec (\frac{v}{2}\
) (53760 (2 \cos (v)+2 \cos (2 v)+2 \cos (3 v)\\
+2 \cos (5 \ v)+1) \sin ^3(\frac{v}{2})+v ((304823 \cos \
(\frac{v}{2})-828603 \cos (\frac{3 v}{2})\\
+554639 \cos (\frac{5 v}{2})-272779 \cos \
(\frac{7 v}{2})) v^2+6720 (7 \cos (\frac{5 \
v}{2})\\
-15 \cos (\frac{7 v}{2})+2 (9 \cos (\
\frac{9 v}{2})-6 \cos (\frac{11 v}{2})+\cos \
(\frac{13 v}{2})))))
\nonumber \\
b_{2,num}^1=\csc ^7(\frac{v}{2}) (v ((393310 \
\cos (2 v)-211395 \cos (3 v)\\
+272779 \cos (4 v)) v^2+(122851 \
v^2+1088640) \cos (v)\\
-15 (-1817 v^2+74368 \cos (2 v)-65856 \
\cos (3 v)
+45248 \cos (4 v)\\
-19040 \cos (5 v)+1344 \cos (6 v)+224 \cos \
(7 v)+36288)) \csc ^3(\frac{v}{2})\\
-26880 (10 \
\cos (v)+10 \cos (2 v)+9 \cos (3 v)+2 \cos (4 v)\\
+8 \cos (5 v)+\cos (6 \
v)+5) \sec (\frac{v}{2}))
\nonumber \\
b_{3,num}^1=\csc ^7(\frac{v}{2}) (26880 (90 \cos \
(v)+90 \cos (2 v)+74 \cos (3 v)\\
+32 \cos (4 v)+58 \cos (5 v)+16 \cos \
(6 v)+45) \sec (\frac{v}{2})\\
-v ((2519110 \cos (2 \
v)-654530 \cos (3 v)+272779 (6 \cos (4 v)\\
+\cos (5 v))) \
v^2+(1288471 v^2+4717440) \cos (v)+8 (47587 \
v^2\\
-643440 \cos (2 v)+593880 \cos (3 v)-416640 \cos (4 v)\\
+159600 \cos \
(5 v)+18480 \cos (6 v)\\
-6720 \cos (7 v)-294840)) \csc \
^3(\frac{v}{2}))
\nonumber \\
b_{4,num}^1=\csc ^7(\frac{v}{2}) (v (1578002 \
v^2+2 (161489 v^2+725760) \cos (v)\\
+160 (20521 \
v^2-13776) \cos (2 v)+5 (478464-4069 v^2) \cos (3 \ v)\\
+2 (738791 v^2-900480) \cos (4 v)+3 (205341 \
v^2+156800) \cos (5 v)\\
+510720 \cos (6 v)-94080 \cos (7 \
v)-725760) \csc ^3(\frac{v}{2})\\
-107520 (30 \cos \
(v)+30 \cos (2 v)+23 \cos (3 v)+14 \cos (4 v)\\
+16 \cos (5 v)+7 \cos (6 \
v)+15) \sec (\frac{v}{2}))
\nonumber \\
b_{5,num}^1=-\csc ^7(\frac{v}{2}) (v (923342 \
v^2+28 (343339 v^2-181440) \cos (v)\\
+20 (240593 \ v^2+103488) \cos (2 v)+5 (832859 v^2+18816) \cos (3 \ v)\\
+2 (1574159 v^2-376320) \cos (4 v)+3 (913431 \
v^2-313600) \cos (5 v)\\
+2446080 \cos (6 v)-376320 \cos (7 \
v)+2540160) \csc ^3(\frac{v}{2})\\
-376320 (30 \cos \
(v)+30 \cos (2 v)+22 \cos (3 v)+16 \cos (4 v)\\
+14 \cos (5 v)+8 \cos (6 \
v)+15) \sec (\frac{v}{2}))
\nonumber \\
b_{6,num}^1=\csc ^7(\frac{v}{2}) (v (7 \
(65953 v^2+181440)+7 (379709 v^2\\
-362880) \cos \ (v)+2 (694753 v^2
+799680) \cos (2 v)+4 (360653 \ v^2\\
-199920) \cos (3 v)+15 (52889 v^2+21952) \cos (4 \ v)\\
+(875393 v^2-540960) \cos (5 v)+799680 \cos (6 v)\\
-117600 \
\cos (7 v)) \csc ^3(\frac{v}{2})-188160 (18 \cos \
(v)+18 \cos (2 v)\\
+13 \cos (3 v)+10 \cos (4 v)+8 \cos (5 v)+5 \cos (6 \
v)+9) \sec (\frac{v}{2}))
\nonumber \\
b_{1,denom}^1=6881280 v^3
\nonumber \\
b_{2,denom}^1=1720320 v^3
\nonumber \\
b_{3,denom}^1=3440640 v^3
\nonumber \\
b_{4,denom}^1=1720320 v^3
\nonumber \\
b_{5,denom}^1=3440640 v^3
\nonumber \\
b_{6,denom}^1=860160 v^3
\nonumber
\ec

\bc
\ds b^1_{T,1}=\frac{90987349}{53222400}-\frac{16301796103 \
v^2}{145297152000}+\frac{1532031563 \
v^4}{268240896000}\\
\ds -\frac{31987133939 \
v^6}{592812380160000}+\frac{5466168990203 \
v^8}{2838385676206080000}+\ldots
\nonumber \\
\ds b^1_{T,2}=-\frac{114798419}{26611200}+\frac{16301796103 \
v^2}{14529715200}-\frac{28198975459 \
v^4}{249080832000}\\
\ds +\frac{88492028011 \
v^6}{11856247603200}-\frac{533071354889581 v^8}{1419192838103040000}+\ldots
\nonumber \\
\ds b^1_{T,3}=\frac{270875723}{17740800}-\frac{16301796103 \
v^2}{3228825600}+\frac{117200518583 \
v^4}{166053888000}\\
\ds -\frac{2285150767769 \
v^6}{39520825344000}+\frac{308808361613933 \
v^8}{105125395415040000}-\ldots
\nonumber \\
\ds b^1_{T,4}=-\frac{67855831}{2217600}+\frac{16301796103 \
v^2}{1210809600}-\frac{36423021893 \
v^4}{16144128000}\\
\ds +\frac{989757481543 \
v^6}{4940103168000}-\frac{1207392034807669 \
v^8}{118266069841920000}\ldots
\nonumber \\
\ds b^1_{T,5}=\frac{50277247}{985600}-\frac{16301796103 \
v^2}{691891200}+\frac{360410789243 \
v^4}{83026944000}\\
\ds -\frac{25035151487 \
v^6}{62731468800}+\frac{125088382253381 \
v^8}{6143691939840000}-\frac{2919882783787129 \
v^{10}}{4460320348323840000}\ldots
\nonumber \\
\ds b^1_{T,6}=-\frac{253491379}{4435200}+\frac{16301796103 \
v^2}{576576000}-\frac{222767191987 \
v^4}{41513472000}\\
\ds +\frac{3516569786117 \
v^6}{7057290240000}-\frac{95477706318691 \
v^8}{3754478407680000}+\ldots
\nonumber \\
\ec

Method \textit{PF-D2}:
\bc
b_{1,num}^2=-30720 (v^2-1) \cos (v)-10 (-387 \csc \
^6(\frac{v}{2})\\
+72 (128 \cos (v)+127) \csc ^2(v)+2048 \
(v^2-3))-(8 (-47 v^4\\
-2040 v^2-(40 \
((22 \cos (v)+22 \cos (2 v)+9 \cos (3 v)+14 \cos (4 v)\\
+8 \cos (5 \
v)-5 \cos (6 v)+11) \tan ^3(\frac{v}{2})-27 v^3) v)/(\
\cos (v)-1)\\
+(4817 v^4+3120 v^2\\
+4410) \cos (v)-2205 \cos (2 \
v)-2205))/((\cos (v)-1)^4)
\nonumber \\
b_{2,num}^2=-1350 v^4 \csc ^{10}(\frac{v}{2})+225 v^2 \
(59 v^2+12) \csc ^8(\frac{v}{2})\\
-5 (7361 \
v^4+3696 v^2-540) \csc ^6(\frac{v}{2})+6 (4817 \
v^4+3120 v^2\\
-7710) \csc ^4(\frac{v}{2})+264780 \csc \
^2(\frac{v}{2})+30720 (3 v^2-7) \cos (v)\\
+15360 \
(v^2-2) \cos (2 v)+30 v \cot (\frac{v}{2}) \
(\frac{1}{8} (-2436 \cos (v)\\
+1921 \cos (2 v)+1595) \csc ^6(\
\frac{v}{2})+8322)+10 (5632 v^2\\
-1024 (8 \cos (v)+23) \
\sin (v) v\\
\ds +\frac{6 (65 \cos (v)+66) \tan (\frac{v}{2}) v}{\
\cos (v)+1}-\frac{36}{\cos (v)+1}-39936)
\nonumber \\
b_{3,num}^2=12150 v^4 \csc ^{10}(\frac{v}{2})-675 v^2 \
(191 v^2+36) \csc ^8(\frac{v}{2})\\
+15 \
(28499 v^4+12432 v^2-1620) \csc ^6(\frac{v}{2})-36450 v \cot (\frac{v}{2}) \csc ^6(\frac{v}{2})\\
-12 (44896 v^4+26400 v^2-35955) \csc \
^4(\frac{v}{2})\\
+381960 v \cot (\frac{v}{2}) \
\csc ^4(\frac{v}{2})+12 (19268 v^4+12480 v^2\\
-219885) \csc ^2(\frac{v}{2})-817950 v \cot (\frac{v}{2}\
) \csc ^2(\frac{v}{2})-540 \sec ^2(\frac{v}{2})\\
-5120 (73 v^2-831)-61440 (9 v^2-40) \cos \
(v)\\
-61440 (3 v^2-8) \cos (2 v)-10240 (v^2-3) \
\cos (3 v)\\
-1512900 v \cot (\frac{v}{2})+1884160 v \sin \
(v)+532480 v \sin (2 v)\\
+30720 v \sin (3 v)+90 v \sec \
^2(\frac{v}{2}) \tan (\frac{v}{2})+7740 v \tan \
(\frac{v}{2})
\nonumber \\
b_{4,num}^2=-4050 v^4 \csc ^{10}(\frac{v}{2})+675 v^2 \
(67 v^2+12) \csc ^8(\frac{v}{2})\\
-15 (11113 \
v^4+4464 v^2-540) \csc ^6(\frac{v}{2})+12150 v \cot \
(\frac{v}{2}) \csc ^6(\frac{v}{2})\\
+6 \
(43259 v^4+24000 v^2-24570) \csc ^4(\frac{v}{2})-133920 v \cot (\frac{v}{2}) \csc ^4(\frac{v}{2}\
)\\
+347250 v \cot (\frac{v}{2}) \csc \
^2(\frac{v}{2})+8 (4817 v^4+15840 v^2-197280)\\
+7680 (19 v^2-125) \cos (v)+30720 (2 v^2-7) \cos (2 v)\\
+7680 (v^2-3) \cos (3 v)+292500 v \cot \
(\frac{v}{2})-555520 v \sin (v)\\
-194560 v \sin (2 v)-23040 \
v \sin (3 v)-30 v \sec ^2(\frac{v}{2}) \tan \
(\frac{v}{2})\\
\ds -3180 v \tan (\frac{v}{2})+\frac{72 (4817 v^4+3120 v^2-26165)}{\cos \
(v)-1}+\frac{360}{\cos (v)+1}
\nonumber \\
b_{5,num}^2=28350 v^4 \csc ^{10}(\frac{v}{2})-14175 v^2 \
(23 v^2+4) \csc ^8(\frac{v}{2})\\
+45 (28269 \
v^4+10864 v^2-1260) \csc ^6(\frac{v}{2})-85050 v \
\cot (\frac{v}{2}) \csc ^6(\frac{v}{2})\\
-12 \
(182576 v^4+98160 v^2-87255) \csc ^4(\frac{v}{2})\\
+965160 v \cot (\frac{v}{2}) \csc ^4(\frac{v}{2}\
)+180 (9634 v^4+6240 v^2\\
-38111) \csc \
^2(\frac{v}{2})-2753430 v \cot (\frac{v}{2}) \
\csc ^2(\frac{v}{2})-256 (2011 v^4\\
+3870 v^2-45810)-15360 (61 v^2-479) \cos (v)-30720 (13 v^2\\
-56) \cos (2 v)-76800 (v^2-3) \cos (3 v)-1054260 v \cot \
(\frac{v}{2})\\
+3537920 v \sin (v)+1372160 v \sin (2 \
v)+230400 v \sin (3 v)\\
\ds -30 v \sec ^2(\frac{v}{2}) \tan \
(\frac{v}{2})-2100 v \tan (\frac{v}{2})+\frac{360}{\cos (v)+1}
\nonumber
\ec

\bc
b_{6,num}^2=93048 v^4+156160 v^2+1042425 \csc ^2(\frac{v}{2})\\
+15360 (9 v^2-74) \cos (v)-268800 \cos (2 v)-38400 \
\cos (3 v)\\
-(9 (7560 \csc ^2(\frac{v}{2})+14963) v^4+187840 v^2-4 (237663 v^4+77080 v^2\\
+298620) \cos \
(v)+5 ((34095 v^4-4000 v^2+63126) \cos (2 v)\\
-2 v \
(v (14451 v^2-5200) \cos (3 v)+8 (5 v (140 \cos (4 \
v)-16 \cos (5 v)\\
-7 \cos (6 v)+2 \cos (7 v))-666 \sin (v)+600 \sin (2 \
v)-298 \sin (3 v)\\
+238 \sin (4 v)+40 \sin (5 v)+80 \sin (6 v)-30 \sin \
(7 v))\\
+36 (\sec ^2(\frac{v}{2})+94) \tan (\
\frac{v}{2})))\\
\ds +878850)/((\cos \
(v)-1)^4)-\frac{270}{\cos (v)+1}-1793280
\nonumber\\
b_{1,denom}^2=5120 v^4
\nonumber \\
b_{2,denom}^2=2560 v^4
\nonumber \\
b_{3,denom}^2=5120 v^4
\nonumber \\
b_{4,denom}^2=640 v^4
\nonumber \\
b_{5,denom}^2=2560 v^4
\nonumber \\
b_{6,denom}^2=320 v^4
\nonumber
\ec

\bc
\ds b^2_{T,1}=\frac{90987349}{53222400}-\frac{16301796103 \
v^2}{96864768000}+\frac{1568734969 \
v^4}{232475443200}\\
\ds -\frac{112423833619 \
v^6}{889218570240000}-\frac{128168340031 \
v^8}{189225711747072000}-\ldots
\nonumber \\
\ds b^2_{T,2}=-\frac{114798419}{26611200}+\frac{16301796103 \
v^2}{9686476800}-\frac{10540703911 \
v^4}{44706816000}\\
\ds +\frac{871935134531 \
v^6}{44460928512000}-\frac{463103326062011 \
v^8}{473064279367680000}+\ldots
\nonumber \\
\ds b^2_{T,3}=\frac{270875723}{17740800}-\frac{16301796103 \
v^2}{2152550400}+\frac{639312597971 \
v^4}{387459072000}\\
\ds -\frac{2472777112313 \
v^6}{11856247603200}+\frac{58110459403753 \
v^8}{3185618042880000}-\ldots
\nonumber \\
\ds b^2_{T,4}=-\frac{67855831}{2217600}+\frac{16301796103 \
v^2}{807206400}-\frac{89147839889 \
v^4}{16144128000}\\
\ds +\frac{1282661397349 \
v^6}{1482030950400}-\frac{3541527647083319 \
v^8}{39422023280640000}+\ldots
\nonumber \\
\ds b^2_{T,5}=\frac{50277247}{985600}-\frac{16301796103 \
v^2}{461260800}+\frac{300047111873 \
v^4}{27675648000}\\
\ds -\frac{18727460555953 \
v^6}{9880206336000}+\frac{4749826752625121 \
v^8}{22526870446080000}-\ldots
\nonumber \\
\ds b^2_{T,6}=-\frac{253491379}{4435200}+\frac{16301796103 \
v^2}{384384000}-\frac{37309797113 \
v^4}{2767564800}\\
\ds +\frac{90333884682737 \
v^6}{37050773760000}-\frac{69221174089601 \
v^8}{250298560512000}+\ldots
\nonumber
\ec

Method \textit{PF-D3}:
\bc
b_{1,num}^3=-108 v^5 \csc ^{10}(\frac{v}{2})+9 v^3 \
(53 v^2+12) \csc ^8(\frac{v}{2})\\
+486 v \csc \
^6(\frac{v}{2})-2484 v \csc ^4(\frac{v}{2})-6144 v (4 v^2-9) \cos (v)\\
+3 \cot (\frac{v}{2}) (99 v^2 \csc ^6(\frac{v}{2})+9 (16-9 v^2) \csc ^4(\frac{v}{2})-8 (113 v^2\\
+276) \csc \
^2(\frac{v}{2})-11266 v^2+10032)-6144 (v^3-6 \
v+4 \sin (v))\\
+2 v^2 (26624 \sin (v)-5635 \tan \
(\frac{v}{2}))+5520 \tan (\frac{v}{2})\\
-9 \
v \sec ^4(\frac{v}{2}) (v \tan (\frac{v}{2})-6)-6 \sec ^2(\frac{v}{2}) ((61 \
v^2+24) \tan (\frac{v}{2})\\
\ds -363 v)+\frac{32508 \
v}{\cos (v)-1}
\nonumber \\
b_{2,num}^3=540 v^5 \csc ^{10}(\frac{v}{2})-45 v^3 \
(71 v^2+12) \csc ^8(\frac{v}{2})\\
+18 v (212 \
v^4+48 v^2-135) \csc ^6(\frac{v}{2})-1485 v^2 \cot \
(\frac{v}{2}) \csc ^6(\frac{v}{2})\\
+15984 v \csc \
^4(\frac{v}{2})+27 (131 v^2-80) \cot \
(\frac{v}{2}) \csc ^4(\frac{v}{2})\\
+55782 v \csc \
^2(\frac{v}{2})+6 (1777 v^2+5952) \cot \
(\frac{v}{2}) \csc ^2(\frac{v}{2})\\
+54 v \sec ^4\
(\frac{v}{2})+6144 v (5 v^2-39)+6144 v (8 \
v^2-33) \cos (v)\\
+18432 v (v^2-4) \cos (2 v)+6 \
(20035 v^2-31704) \cot (\frac{v}{2})\\
-20480 \
(7 v^2-6) \sin (v)+3072 (12-19 v^2) \sin (2 \
v)\\
-9 v^2 \sec ^4(\frac{v}{2}) \tan (\frac{v}{2})-12 (23 v^2+12) \sec ^2(\frac{v}{2}) \tan \
(\frac{v}{2})\\
\ds +14 (456-337 v^2) \tan \
(\frac{v}{2})+\frac{3276 v}{\cos (v)+1}
\nonumber \\
b_{3,num}^3=-1620 v^5 \csc ^{10}(\frac{v}{2})+135 v^3 \
(85 v^2+12) \csc ^8(\frac{v}{2})\\
-6 v (3872 \
v^4+768 v^2-1215) \csc ^6(\frac{v}{2})+4455 v^2 \cot \
(\frac{v}{2}) \csc ^6(\frac{v}{2})\\
+36 v \
(424 v^4+96 v^2-1563) \csc ^4(\frac{v}{2})+9 \
(720\\
-1781 v^2) \cot (\frac{v}{2}) \csc ^4(\
\frac{v}{2})-96642 v \csc ^2(\frac{v}{2})-24 \
(685 v^2\\
+4716) \cot (\frac{v}{2}) \csc ^2(\
\frac{v}{2})-54 v \sec ^4(\frac{v}{2})+1024 v \
(606-67 v^2)\\
-49152 v (2 v^2-13) \cos (v)-49152 \
v (v^2-5) \cos (2 v)\\
+2048 v (21-4 v^2) \cos (3 \
v)+6 (111720-46199 v^2) \cot (\frac{v}{2})\\
+32768 (11 v^2-15) \sin (v)+4096 (41 v^2-36) \sin (2 v)\\
+4096 (7 v^2-6) \sin (3 v)+9 v^2 \sec \
^4(\frac{v}{2}) \tan (\frac{v}{2})+6 (61 \
v^2\\
\ds +24) \sec ^2(\frac{v}{2}) \tan (\frac{v}{2})+10 (1055 v^2-552) \tan (\frac{v}{2})-\frac{4356 v}{\cos (v)+1}
\nonumber \\
b_{4,num}^3=1620 v^5 \csc ^{10}(\frac{v}{2})-135 v^3 \
(95 v^2+12) \csc ^8(\frac{v}{2})\\
+18 v \
(1844 v^4+336 v^2-405) \csc ^6(\frac{v}{2})-4455 v^2 \cot (\frac{v}{2}) \csc ^6(\frac{v}{2}\
)\\
-216 v (169 v^4+36 v^2-288) \csc \
^4(\frac{v}{2})+27 (737 v^2\\
-240) \cot \
(\frac{v}{2}) \csc ^4(\frac{v}{2})+18 v \
(848 v^4+192 v^2+2229) \csc ^2(\frac{v}{2})\\
+18 \
(35 v^2+6528) \cot (\frac{v}{2}) \csc \
^2(\frac{v}{2})-54 v \sec ^4(\frac{v}{2})+384 v \
(133 v^2\\
-1350)+768 v (112 v^2-843) \cos \
(v)+3072 v (13 v^2-81) \cos (2 v)\\
+1536 v (8 v^2-45) \cos (3 v)+768 v (v^2-6) \cos (4 v)+54 (4447 \
v^2\\
-13800) \cot (\frac{v}{2})+512 (1146-659 v^2\
) \sin (v)-512 (289 v^2\\
-348) \sin (2 v)+1024 \
(42-43 v^2) \sin (3 v)+256 (12\\
\ds -11 v^2) \sin (4 \
v)+(9 v^2 \sec ^4(\frac{v}{2})+5258 v^2+\frac{552 \
v^2+288}{\cos (v)+1}\\
\ds -6384) \tan (\frac{v}{2})-\frac{3276 v}{\cos (v)+1}
\nonumber
\ec

\bc
b_{5,num}^3=-11340 v^5 \csc ^{10}(\frac{v}{2})+945 v^3 \
(101 v^2+12) \csc ^8(\frac{v}{2})\\
-18 v \
(15472 v^4+2688 v^2-2835) \csc ^6(\frac{v}{2})\\
+31185 v^2 \cot (\frac{v}{2}) \csc \
^6(\frac{v}{2})+108 v (3500 v^4+720 v^2-4263) \
\csc ^4(\frac{v}{2})\\
+189 (240-823 v^2) \cot \
(\frac{v}{2}) \csc ^4(\frac{v}{2})+864 (87 \
v^2\\
-973) \cot (\frac{v}{2}) \csc ^2(\frac{v}{2}\
)+54 v \sec ^4(\frac{v}{2})+384 v (159 v^4-796 \
v^2\\
+8568)+3072 v (1395-172 v^2) \cos (v)-24576 v \
(11 v^2\\
-75) \cos (2 v)+3072 v (171-28 v^2) \cos \
(3 v)\\
-12288 v (v^2-6) \cos (4 v)+18 (307848-88943 v^2\
) \cot (\frac{v}{2})\\
+2048 (1055 v^2-2154) \
\sin (v)+2048 (505 v^2-708) \sin (2 v)\\
+2048 (155 \
v^2-174) \sin (3 v)+4096 (11 v^2-12) \sin (4 v)\\
-9 \
v^2 \sec ^4(\frac{v}{2}) \tan (\frac{v}{2})-6 \
(61 v^2+24) \sec ^2(\frac{v}{2}) \tan \
(\frac{v}{2})+10 (552\\
\ds -1019 v^2) \tan \
(\frac{v}{2})+\frac{36 v (256 (53 v^2+12) \
v^2+1479)}{\cos (v)-1}+\frac{4356 v}{\cos (v)+1}
\nonumber\\
b_{6,num}^3=2268 v^5 \csc ^{10}(\frac{v}{2})-189 v^3 \
(103 v^2+12) \csc ^8(\frac{v}{2})\\
+6 v \
(9820 v^4+1680 v^2-1701) \csc ^6(\frac{v}{2})-6237 v^2 \cot (\frac{v}{2}) \csc ^6(\frac{v}{2}\
)\\
-144 v (590 v^4+120 v^2-651) \csc \
^4(\frac{v}{2})+63 (511 v^2\\
-144) \cot \
(\frac{v}{2}) \csc ^4(\frac{v}{2})+18 v \
(3392 v^4+768 v^2-677) \csc ^2(\frac{v}{2})\\
+6 \
(28224-3463 v^2) \cot (\frac{v}{2}) \csc \
^2(\frac{v}{2})+54 v \sec ^4(\frac{v}{2})\\
-128 v \
(141 v^4-460 v^2+5016)+3072 v (32 v^2-269) \cos \
(v)\\
+6144 v (9 v^2-62) \cos (2 v)+2048 v (8 v^2-51) \cos (3 v)\\
+3072 v (v^2-6) \cos (4 v)+1002 (317 \
v^2-1128) \cot (\frac{v}{2})\\
-6144 (67 v^2-146) \sin (v)-3072 (69 v^2-100) \sin (2 v)\\
-12288 (5 \
v^2-6) \sin (3 v)+1024 (12-11 v^2) \sin (4 \
v)\\
\ds -(9 v^2 \sec ^4(\frac{v}{2})+5438 v^2+\frac{552 \
v^2+288}{\cos (v)+1}-6384) \tan (\frac{v}{2})\\
\ds +\frac{3276 v}{\cos (v)+1}
\nonumber \\
b_{1,denom}^3=3072 v^5
\nonumber \\
b_{2,denom}^3=1536 v^5
\nonumber \\
b_{3,denom}^3=1024 v^5
\nonumber \\
b_{4,denom}^3=384 v^5
\nonumber \\
b_{5,denom}^3=1536 v^5
\nonumber \\
b_{6,denom}^3=256 v^5
\nonumber\ec

\bc
\ds b^3_{T,1}=\frac{90987349}{53222400}-\frac{16301796103 \
v^2}{72648576000}+\frac{1273143229 v^4}{193729536000}\\
\ds -\frac{421960559 \
v^6}{1221454080000}-\frac{124506164597657 \
v^8}{5676771352412160000}-\ldots
\nonumber \\
\ds b^3_{T,2}=-\frac{114798419}{26611200}+\frac{16301796103 \
v^2}{7264857600}-\frac{38969308351 v^4}{96864768000}\\
\ds +\frac{5604849953 \
v^6}{158789030400}-\frac{5242236508523657 \
v^8}{2838385676206080000}+\ldots
\nonumber \\
\ds b^3_{T,3}=\frac{270875723}{17740800}-\frac{16301796103 \
v^2}{1614412800}+\frac{578948920601 \
v^4}{193729536000}\\
\ds -\frac{733132845253 \
v^6}{1482030950400}+\frac{3064860132866873 \
v^8}{57341124771840000}-\ldots
\nonumber \\
\ds b^3_{T,4}=-\frac{67855831}{2217600}+\frac{16301796103 \
v^2}{605404800}-\frac{247322293877 \
v^4}{24216192000}\\
\ds +\frac{844539338551 \
v^6}{370507737600}-\frac{80127849554001773 \
v^8}{236532139683840000}+\ldots
\nonumber \\
\ds b^3_{T,5}=\frac{50277247}{985600}-\frac{16301796103 \
v^2}{345945600}+\frac{279925886083 \
v^4}{13837824000}\\
\ds -\frac{644885438797 \
v^6}{123502579200}+\frac{122198747647321067 v^8}{135161222676480000}-\ldots
\nonumber \\
\ds b^3_{T,6}=-\frac{253491379}{4435200}+\frac{16301796103 \
v^2}{288288000}-\frac{1220113637 v^4}{48384000}\\
\ds +\frac{55089606671 \
v^6}{8096760000}-\frac{1634963205351829 \
v^8}{1325110026240000}+\ldots
\nonumber \\
\ec

Method \textit{PF-D4}:
\bc
b_{1,num}^4=216 v^6 \csc ^{10}(\frac{v}{2})-810 v^2 \csc \
^6(\frac{v}{2})-54 (157 v^2\\
+80) \csc \
^4(\frac{v}{2})+36 (1840-2257 v^2) \csc \
^2(\frac{v}{2})-12288 (20 v^4\\
-71 v^2+20) \cos \
(v)-3 v \cot (\frac{v}{2}) ((4682 v^2+9 \csc \
^2(\frac{v}{2}) (5 \csc ^2(\frac{v}{2}) \
v^2\\
+54 v^2+96)+5856) \csc ^2(\frac{v}{2})+2520 \
(31 v^2-48))\\
-45 v^2 \sec ^6(\frac{v}{2}) \
(v \tan (\frac{v}{2})-6)+56 v (1024 \
(11 v^2-12) \sin (v)\\
+135 (48-31 v^2) \tan \
(\frac{v}{2}))-18 \sec ^4(\frac{v}{2}) \
(-373 v^2+3 (21 v^2\\
+16) \tan (\frac{v}{2}) \
v-80)-2 \sec ^2(\frac{v}{2}) (-39246 \
v^2\\
+(6733 v^2+6384) \tan (\frac{v}{2}) v+28320)
\nonumber \\
b_{2,num}^4=-1080 v^6 \csc ^{10}(\frac{v}{2})+2160 v^6 \
\csc ^8(\frac{v}{2})+4050 v^2 \csc ^6(\frac{v}{2})\\
+54 (617 v^2+400) \csc ^4(\frac{v}{2})+293508 v^2 \csc ^2(\frac{v}{2})+12288 (20 \
v^4\\
-71 v^2+20) \cos (v)+24576 (10 v^4-59 v^2+20) \
\cos (2 v)\\
+3 v \cot (\frac{v}{2}) ((16978 v^2+9 \
\csc ^2(\frac{v}{2}) (25 \csc ^2(\frac{v}{2}) v^2+214 v^2+480)\\
+20064) \csc ^2(\frac{v}{2})+120 (2341 v^2-5136))+(72 (4096 \
v^4\\
-31783 v^2-80 (257 \cos (v)+251) \csc ^2(v)+1024 (4 v^4-35 \
v^2\\
+20) \cos (v)+25600))/(\cos (v)+1)-v (135 v \
(v \tan (\frac{v}{2})\\
-6) \sec \
^6(\frac{v}{2})+18 ((161 v^2+144) \tan \
(\frac{v}{2})-951 v) \sec ^4(\frac{v}{2})\\
+18 (1407 v^2+1616) \tan (\frac{v}{2}) \
\sec ^2(\frac{v}{2})+8 (7168 (11 v^2-12) \
\sin (v)\\
+1024 (107 v^2-156) \sin (2 v)+15 (3217 \
v^2-9168) \tan (\frac{v}{2})))
\nonumber \\
b_{3,num}^4=9720 v^6 \csc ^{10}(\frac{v}{2})-34560 v^6 \
\csc ^8(\frac{v}{2})+270 v^2 (128 v^4\\
-135) \csc \
^6(\frac{v}{2})-6075 v^3 \cot (\frac{v}{2}) \
\csc ^6(\frac{v}{2})-162 (1459 v^2\\
+1200) \csc \
^4(\frac{v}{2})-54 v (767 v^2+2160) \cot (\
\frac{v}{2}) \csc ^4(\frac{v}{2})\\
+108 \
(32720-18831 v^2) \csc ^2(\frac{v}{2})-234 v \
(1511 v^2\\
+1488) \cot (\frac{v}{2}) \csc \
^2(\frac{v}{2})+3510 v^2 \sec ^6(\frac{v}{2})+90 (701 v^2\\
+208) \sec ^4(\frac{v}{2})-12288 (98 v^4-825 v^2+540)-12288 (220 v^4\\
-1717 \
v^2+940) \cos (v)-98304 (10 v^4-59 v^2+20) \cos (2 \
v)\\
-24576 (20 v^4-147 v^2+60) \cos (3 v)+360 v \
(45744\\
-16567 v^2) \cot (\frac{v}{2})+8192 v \
(1279 v^2-2868) \sin (v)+32768 v (107 v^2\\
-156) \
\sin (2 v)+147456 v (13 v^2-24) \sin (3 v)-v (585 v^2 \
\sec ^6(\frac{v}{2})\\
+(42991 v^2+(32281 v^2+46128\
) \cos (v)+57360) \sec ^4(\frac{v}{2})\\
\ds +120 \
(8741 v^2-38928)) \tan (\frac{v}{2})+\frac{24 (30457 v^2-69040)}{\cos (v)+1}
\nonumber \\
b_{4,num}^4=-3240 v^6 \csc ^{10}(\frac{v}{2})+15120 v^6 \
\csc ^8(\frac{v}{2})-810 v^2 (32 v^4\\
-15) \csc \
^6(\frac{v}{2})+2025 v^3 \cot (\frac{v}{2}) \
\csc ^6(\frac{v}{2})+54 (320 v^6+1179 v^2\\
+1200) \
\csc ^4(\frac{v}{2})+54 v (209 v^2+720) \cot \
(\frac{v}{2}) \csc ^4(\frac{v}{2})\\
+108 \
(5237 v^2-11440) \csc ^2(\frac{v}{2})+18 v \
(5461 v^2\\
+3888) \cot (\frac{v}{2}) \csc \
^2(\frac{v}{2})+270 v^2 \sec ^6(\frac{v}{2})+90 \
(41 v^2+16) \sec ^4(\frac{v}{2})\\
+76800 (5 \
v^4-54 v^2+48)+6144 (100 v^4-819 v^2+480) \cos \
(v)\\
+73728 (5 v^4-39 v^2+20) \cos (2 v)+6144 (20 \
v^4-147 v^2\\
+60) \cos (3 v)+6144 (5 v^4-41 v^2+20) \
\cos (4 v)+360 v (4705 v^2\\
-14928) \cot (\frac{v}{2})-36864 v (67 v^2-154) \sin (v)-86016 v (17 \
v^2\\
-36) \sin (2 v)-36864 v (13 v^2-24) \sin (3 \
v)-2048 v (61 v^2\\
-132) \sin (4 v)-v \
\ds (((1087+\frac{90}{\cos (v)+1}) v^2+(457 \
v^2\\
\ds +1776) \cos (v)+2640) \sec ^4(\frac{v}{2})+120 (349 v^2\\
\ds -2832)) \tan (\frac{v}{2})+\frac{24 (349 v^2-5680)}{\cos (v)+1}
\nonumber
\ec

\bc
b_{5,num}^4=22680 v^6 \csc ^{10}(\frac{v}{2})-120960 v^6 \
\csc ^8(\frac{v}{2})+4050 v^2 (64 v^4\\
-21) \csc \
^6(\frac{v}{2})-14175 v^3 \cot (\frac{v}{2}) \
\csc ^6(\frac{v}{2})-54 (5120 v^6+7077 v^2\\
+8400) \csc ^4(\frac{v}{2})-378 v (181 v^2+720) \
\cot (\frac{v}{2}) \csc ^4(\frac{v}{2})+108 \
(1280 v^6\\
-33131 v^2+82320) \csc ^2(\frac{v}{2})-378 v (1643 v^2+784) \cot (\frac{v}{2}) \
\csc ^2(\frac{v}{2})\\
-1890 v^2 \sec ^6(\frac{v}{2})-126 (277 v^2+80) \sec ^4(\frac{v}{2})-15360 (160 v^4\\
-1833 v^2+1740)-30720 (120 \
v^4-1079 v^2+740) \cos (v)\\
-6144 (400 v^4-3271 v^2+1780) \cos (2 v)-6144 (140 v^4-1113 v^2\\
+540) \cos (3 \
v)-49152 (5 v^4-41 v^2+20) \cos (4 v)-6144 (4 v^4\\
-35 \
v^2+20) \cos (5 v)+360 v (102576-30047 v^2) \cot \
(\frac{v}{2})\\
+20480 v (755 v^2-1956) \sin \
(v)+4096 v (2429 v^2-5412) \sin (2 v)\\
+12288 v (281 \
v^2-588) \sin (3 v)+16384 v (61 v^2-132) \sin (4 \
v)\\
+20480 v (5 v^2-12) \sin (5 v)+v \
(((25429+\frac{630}{\cos (v)+1}) v^2+(19507 \
v^2\\
\ds +26256) \cos (v)+32304) \sec ^4(\frac{v}{2})+120 (5347 v^2-20976)) \tan (\frac{v}{2})\\
\ds +\frac{885120-443976 v^2}{\cos (v)+1}
\nonumber \\
b_{6,num}^4=-13608 v^6 \csc ^{10}(\frac{v}{2})+75600 v^6 \
\csc ^8(\frac{v}{2})+270 v^2 (189\\
-640 v^4) \csc \
^6(\frac{v}{2})+8505 v^3 \cot (\frac{v}{2}) \
\csc ^6(\frac{v}{2})+162 (1280 v^6\\
+1337 v^2+1680) \csc ^4(\frac{v}{2})+378 v (103 v^2+432) \
\cot (\frac{v}{2}) \csc ^4(\frac{v}{2})\\
+126 v \
(2863 v^2+1104) \cot (\frac{v}{2}) \csc \
^2(\frac{v}{2})-1890 v^2 \sec ^6(\frac{v}{2})\\
-126 (253 v^2+80) \sec ^4(\frac{v}{2})+3072 (18 v^6+400 v^4-4575 v^2+4500)\\
+30720 (80 \
v^4-753 v^2+540) \cos (v)+30720 (40 v^4-327 v^2\\
+180) \
\cos (2 v)+6144 (100 v^4-819 v^2+420) \cos (3 v)+24576 \
(5 v^4\\
-41 v^2+20) \cos (4 v)+6144 (4 v^4-35 v^2+20) \cos (5 v)\\
+360 v (17495 v^2-61104) \cot \
(\frac{v}{2})-552960 v (19 v^2-52) \sin \
(v)\\
-552960 v (9 v^2-20) \sin (2 v)-12288 v (203 \
v^2-444) \sin (3 v)\\
-8192 v (61 v^2-132) \sin (4 \
v)-20480 v (5 v^2-12) \sin (5 v)\\
+v (63 (48 \
(v^2+1)+(43 v^2+48) \cos (v)) \sec \
^6(\frac{v}{2})+840 (659 v^2\\
\ds -2928)+\frac{57964 \
v^2+86592}{\cos (v)+1}) \tan (\frac{v}{2})\\
\ds +\frac{216 \
(1280 v^6-19257 v^2+49840)}{\cos \
(v)-1}+\frac{908160-325608 v^2}{\cos (v)+1}
\nonumber \\
b_{1,denom}^4=24576 v^6
\nonumber \\
b_{2,denom}^4=12288 v^6
\nonumber \\
b_{3,denom}^4=24576 v^6
\nonumber \\
b_{4,denom}^4=3072 v^6
\nonumber \\
b_{5,denom}^4=12288 v^6
\nonumber \\
b_{6,denom}^4=6144 v^6
\nonumber
\ec

\bc
\ds b^4_{T,1}=\frac{90987349}{53222400}-\frac{16301796103 \
v^2}{58118860800}+\frac{451826777 \
v^4}{87178291200}\\
\ds -\frac{567827928197 \
v^6}{711374856192000}-\frac{7456559915267 \
v^8}{87334943883264000}-\ldots
\nonumber \\
\ds b^4_{T,2}=-\frac{114798419}{26611200}+\frac{16301796103 \
v^2}{5811886080}-\frac{53423656079 \
v^4}{87178291200}+\\
\ds \frac{3768047264957 \
v^6}{71137485619200}-\frac{1807362562020419 v^8}{567677135241216000}+\ldots
\nonumber \\
\ds b^4_{T,3}=\frac{270875723}{17740800}-\frac{16301796103 \
v^2}{1291530240}+\frac{137191770479 \
v^4}{29059430400}\\
\ds -\frac{6482264498201 \
v^6}{6774998630400}+\frac{177746918811101 \
v^8}{1557413265408000}-\ldots
\nonumber \\
\ds b^4_{T,4}=-\frac{67855831}{2217600}+\frac{16301796103 \
v^2}{484323840}-\frac{4393734833 v^4}{269068800}\\
\ds +\frac{3998343544387 \
v^6}{846874828800}-\frac{42449166474915191 \
v^8}{47306427936768000}+\ldots
\nonumber \\
\ds b^4_{T,5}=\frac{50277247}{985600}-\frac{16301796103 \
v^2}{276756480}+\frac{9638045471 v^4}{296524800}\\
\ds -\frac{87748735805873 \
v^6}{7904165068800}+\frac{70694795981650949 \
v^8}{27032244535296000}-\ldots
\nonumber \\
\ds b^4_{T,6}=-\frac{253491379}{4435200}+\frac{16301796103 \
v^2}{230630400}-\frac{84219941177 \
v^4}{2075673600}\\
\ds +\frac{863717128805537 \
v^6}{59281238016000}-\frac{1830923526668191 \
v^8}{500597121024000}+\ldots
\nonumber
\ec

Method \textit{PF-D5}:
\bc
b_{1,num}^5=-\cos ^5(\frac{v}{2}) \csc ^{13}(v) \sin \
^7(\frac{v}{2}) (60 (-5 v^5+3 v^3+3 \sin (v)\\
-6 \
\sin (3 v)+3 \sin (4 v)+3 \sin (5 v)-6 \sin (6 v)+3 \sin (8 v))\\
+\frac{1}{2} v (15 (200 v^4+283 v^2-100) \cos \
(v)+300 (2 \cos (2 v)+6 \cos (3 v)\\
-10 \cos (4 v)+\cos (5 v)+12 \cos \
(6 v)-2 \cos (7 v)-4 \cos (8 v))\\
+v (v (15 (366 \cos (3 \
v)+686 \cos (4 v)-503 \cos (5 v)-480 \cos (6 v)\\
+142 \cos (7 v)+116 \
\cos (8 v))+2 (15 (72 v^2-155) \cos (2 v)\\
+v (-60 \
v (45 \cos (3 v)+24 \cos (4 v)-25 \cos (5 v)-14 \cos (6 v)\\
+5 \cos (7 \
v)+3 \cos (8 v))+405 \cot (\frac{v}{2})+357 \sin (v)-2530 \
\sin (2 v)\\
+4516 \sin (3 v)+3752 \sin (4 v)-3373 \sin (5 v)-2324 \sin \
(6 v)\\
+770 \sin (7 v)+522 \sin (8 v)+315 \tan (\frac{v}{2})\
)))-60 (49 \sin (v)\\
-44 \sin (2 v)-6 \sin (3 v)+131 \
\sin (4 v)-67 \sin (5 v)-114 \sin (6 v)\\
+28 \sin (7 v)+31 \sin (8 v)))))
\nonumber \\
b_{2,num}^5=\csc ^7(\frac{v}{2}) \sec ^9(\frac{v}{2}) (4 (3875 v^4-519 v^2-90) \cos (v)\\
-36 (147 \
v^4-76 v^2-10 \cos (2 v)+20 \cos (4 v)-20 \cos (5 v)\\
+10 \cos (7 v)-20 \
\cos (8 v)+10 \cos (10 v))+v (v (2 (5117 \
v^2\\
-4254) \cos (2 v)+12 (516 \cos (3 v)+482 \cos (4 v)-158 \cos \
(5 v)\\
+348 \cos (6 v)-317 \cos (7 v)-486 \cos (8 v)+132 \cos (9 v)\\
+137 \
\cos (10 v))+v (3 (-12 \sin (v)-3789 \sin (2 v)+3900 \sin (3 \
v)\\
+2058 \sin (4 v)+708 \sin (5 v)+1632 \sin (6 v)-2142 \sin (7 \
v)\\
-1792 \sin (8 v)+614 \sin (9 v)+461 \sin (10 v))\\
-2 v (-393216 \
v (15 \cos (2 v)+17) \sin ^7(\frac{v}{2}) \cos \
^9(\frac{v}{2})+4476 \cos (3 v)\\
+1918 \cos (4 v)+2288 \cos \
(5 v)+1242 \cos (6 v)-2476 \cos (7 v)\\
-1400 \cos (8 v)+582 \cos (9 \
v)+351 \cos (10 v))))\\
-15360 \cos \
^3(\frac{v}{2}) (186 \cos (v)+156 \cos (2 v)+132 \cos (3 \
v)\\
+78 \cos (4 v)+48 \cos (5 v)+19 \cos (6 v)+98) \sin \
^5(\frac{v}{2})))
\nonumber \\
b_{3,num}^5=-\cos ^6(\frac{v}{2}) \csc ^{14}(v) \sin \
^8(\frac{v}{2}) (7290 \cot (\frac{v}{2}) \
v^4\\
+5670 \tan (\frac{v}{2}) v^4+32 (56 v^2-93) \
\sin (9 v) v^2+12 (250 v^4-751 v^2\\
+100) v+3 (2140 \
v^4+13737 v^2-2220) \cos (v) v-6 (260 v^4\\
-591 v^2+540) \cos (2 v) v+6 (-1702 v^4+2971 v^2+1100) \cos (3 v) \
v\\
-6 (584 v^4-3933 v^2+1020) \cos (4 v) v+15 (236 \
v^4-1395 v^2\\
+228) \cos (5 v) v+12 (182 v^4-1121 v^2+800) \cos (6 v) v+12 (61 v^4\\
-201 v^2-180) \cos (7 v) v+12 \
(10 v^4-171 v^2+60) \cos (8 v) v\\
-60 (8 v^4-53 v^2+20) \cos (9 v) v-24 (10 v^4-93 v^2+90) \cos (10 v) v\\
+6 \
(1327 v^4-974 v^2-60) \sin (v)-4 (1295 v^4+1098 \
v^2\\
-180) \sin (2 v)+8 (3403 v^4+579 v^2-270) \sin (3 \
v)+4 (3256 v^4\\
-4557 v^2+90) \sin (4 v)+2 (-6097 \
v^4+8058 v^2+180) \sin (5 v)\\
-8 (935 v^4-1929 v^2+270) \sin (6 v)-4 (593 v^4+126 v^2\\
-180) \sin (7 v)-4 \
(v^2-3) (199 v^2+30) \sin (8 v)+8 (127 \
v^4-363 v^2\\
+90) \sin (10 v))
\nonumber
\ec

\bc
b_{4,num}^5=5670 v^3 \sec ^8(\frac{v}{2})-945 v^4 \tan \
(\frac{v}{2}) \sec ^8(\frac{v}{2})+270 v \
(257 v^2\\
+160) \sec ^6(\frac{v}{2})-1080 v^2 \
(11 v^2+18) \tan (\frac{v}{2}) \sec \
^6(\frac{v}{2})+4320 v (51 v^2\\
-46) \sec \
^4(\frac{v}{2})-18 (2267 v^4+9816 v^2+2880) \tan (\frac{v}{2}) \sec ^4(\frac{v}{2})\\
+4 \
(-61651 v^4+89664 v^2+570240) \tan (\frac{v}{2}) \sec ^2(\frac{v}{2})\\
+7290 v^3 \csc \
^6(\frac{v}{2})+1215 v^4 \cot (\frac{v}{2}) \
\csc ^6(\frac{v}{2})+78732 v^3 \csc ^4(\frac{v}{2})\\
+81 v^2 (167 v^2+240) \cot (\frac{v}{2}) \
\csc ^4(\frac{v}{2})+756 v (1127 v^2\\
-2640) \csc \
^2(\frac{v}{2})+9 (16297 v^4+2592 v^2-17280) \
\cot (\frac{v}{2}) \csc ^2(\frac{v}{2})\\
+46080 v \
(30 v^4-443 v^2+980)+73728 v (101 v^2-370) \cos \
(v)\\
+36864 v (40 v^4-431 v^2+580) \cos (2 v)+73728 v \
(13 v^2\\
-40) \cos (3 v)+6144 v (30 v^4-307 v^2+340) \cos (4 v)+9 (342973 v^4\\
-1171680 v^2+328320) \cot \
(\frac{v}{2})+73728 (17 v^4-271 v^2\\
+210) \sin \
(v)-98304 (67 v^4-243 v^2+90) \sin (2 v)\\
+16384 (10 \
v^4-147 v^2+90) \sin (3 v)-24576 (33 v^4-107 v^2\\
+30) \
\sin (4 v)+(-3937883 v^4+28223520 v^2-21006720) \tan \
(\frac{v}{2})\\
\ds +\frac{24 v (115673 v^2-893360)}{\cos (v)+1}
\nonumber \\
b_{5,num}^5=39690 v^3 \sec ^8(\frac{v}{2})-6615 v^4 \tan \
(\frac{v}{2}) \sec ^8(\frac{v}{2})+1890 v \
(239 v^2\\
+160) \sec ^6(\frac{v}{2})-1890 v^2 \
(41 v^2+72) \tan (\frac{v}{2}) \sec \
^6(\frac{v}{2})+756 v (1471 v^2\\
-2160) \sec \
^4(\frac{v}{2})-378 (561 v^4+2968 v^2+960) \tan \
(\frac{v}{2}) \sec ^4(\frac{v}{2})\\
+84 v \
(99661 v^2-859120) \sec ^2(\frac{v}{2})-56 \
(26371 v^4-61404 v^2\\
-289440) \tan (\frac{v}{2}) \
\sec ^2(\frac{v}{2})-51030 v^3 \csc ^6(\frac{v}{2})\\
-8505 v^4 \cot (\frac{v}{2}) \csc ^6(\frac{v}{2}\
)-503496 v^3 \csc ^4(\frac{v}{2})-1701 v^2 (51 \
v^2\\
+80) \cot (\frac{v}{2}) \csc ^4(\frac{v}{2})-756 v (7131 v^2-18320) \csc ^2(\frac{v}{2})\\
-189 (4913 v^4+256 v^2-5760) \cot (\frac{v}{2}) \csc ^2(\frac{v}{2})-92160 v (373 v^2\\
-2140)-245760 v (60 v^4-767 v^2+1400) \cos (v)-12288 v \
(1829 v^2\\
-6620) \cos (2 v)-884736 v (5 v^4-54 v^2+70) \cos (3 v)-24576 v (61 v^2\\
-220) \cos (4 v)-49152 v \
(6 v^4-65 v^2+80) \cos (5 v)-9 (2193517 v^4\\
-8002080 \
v^2+2378880) \cot (\frac{v}{2})+245760 (287 \
v^4-1321 v^2\\
+690) \sin (v)-24576 (151 v^4-2446 v^2+1860) \sin (2 v)\\
+16384 (1211 v^4-4317 v^2+1530) \sin (3 \
v)-49152 (5 v^4-82 v^2\\
+60) \sin (4 v)+49152 (27 \
v^4-95 v^2+30) \sin (5 v)+(-24636373 v^4\\
+185682720 \
v^2-157806720) \tan (\frac{v}{2})
\nonumber \\
b_{6,num}^5=119070 v^3 \sec ^8(\frac{v}{2})-19845 v^4 \tan \
(\frac{v}{2}) \sec ^8(\frac{v}{2})\\
+5670 v \
(233 v^2+160) \sec ^6(\frac{v}{2})-45360 v^2 \
(5 v^2+9) \tan (\frac{v}{2}) \sec \
^6(\frac{v}{2})\\
+30240 v (97 v^2-170) \sec \
^4(\frac{v}{2})-126 (4501 v^4+25800 v^2\\
+8640) \tan (\frac{v}{2}) \sec ^4(\frac{v}{2})+28 \
(-151573 v^4+398400 v^2\\
+1745280) \tan (\frac{v}{2}) \sec ^2(\frac{v}{2})+153090 v^3 \csc
^6(\frac{v}{2})\\
+25515 v^4 \cot (\frac{v}{2}) \
\csc ^6(\frac{v}{2})+1462860 v^3 \csc ^4(\frac{v}{2})+2835 v^2 (89 v^2\\
+144) \cot (\frac{v}{2}) \
\csc ^4(\frac{v}{2})+3780 v (4139 v^2-10960) \
\csc ^2(\frac{v}{2})\\
+63 (42779 v^4+480 v^2-51840) \cot (\frac{v}{2}) \csc ^2(\frac{v}{2})+61440 v (400 v^4\\
-6273 v^2+15540)+184320 v (737 \
v^2-3260) \cos (v)\\
+184320 v (150 v^4-1739 v^2+2660) \
\cos (2 v)+368640 v (61 v^2\\
-220) \cos (3 v)+61440 v \
(72 v^4-779 v^2+980) \cos (4 v)\\
+614400 v (v^2-4) \cos (5 v)+30720 v (4 v^4-45 v^2+60) \cos (6 v)\\
+27 \
(2129819 v^4-7933600 v^2+2405760) \cot (\frac{v}{2})+122880 (173 v^4\\
-3222 v^2+3060) \sin (v)-491520 \
(259 v^4-1035 v^2+450) \sin (2 v)\\
+8192 (451 v^4-7365 \
v^2+5490) \sin (3 v)-24576 (809 v^4-2860 v^2
\ec

\bc
+960) \
\sin (4 v)+49152 (2 v^4-35 v^2+30) \sin (5 v)-4096 \
(137 v^4\\
-510 v^2+180) \sin (6 v)+(-71232167 \
v^4+545068320 v^2\\
\ds -484318080) \tan (\frac{v}{2})+\frac{4200 v (11495 v^2-101648)}{\cos (v)+1}
\nonumber\\
b_{1,denom}^5=30 v^7
\nonumber \\
b_{2,denom}^5=393216 v^7
\nonumber \\
b_{3,denom}^5=6 v^7
\nonumber \\
b_{4,denom}^5=6144 v^7
\nonumber \\
b_{5,denom}^5=24576 v^7
\nonumber \\
b_{6,denom}^5=61440 v^7
\nonumber
\ec

\bc
\ds b^5_{T,1}=\frac{90987349}{53222400}-\frac{16301796103 \
v^2}{48432384000}+\frac{3000167803 \
v^4}{1162377216000}-\\
\ds \frac{2795461105681 \
v^6}{1778437140480000}-\frac{423149542577683 \
v^8}{1892257117470720000}-\ldots
\nonumber \\
\ds b^5_{T,2}=-\frac{114798419}{26611200}+\frac{16301796103 \
v^2}{4843238400}-\frac{9164631311 \
v^4}{10567065600}+\\
\ds \frac{1149042063431 \
v^6}{16167610368000}-\frac{1015140295561079 v^8}{189225711747072000}-\ldots
\nonumber \\
\ds b^5_{T,3}=\frac{270875723}{17740800}-\frac{16301796103 \
v^2}{1076275200}+\frac{106131595741 \
v^4}{15498362880}-\\
\ds \frac{193910381341843 \
v^6}{118562476032000}+\frac{12701953007551 \
v^8}{62296530616320}-\ldots
\nonumber \\
\ds b^5_{T,4}=-\frac{67855831}{2217600}+\frac{16301796103 \
v^2}{403603200}-\frac{15415020883 \
v^4}{645765120}+\\
\ds \frac{125531790304181 \
v^6}{14820309504000}-\frac{361081038637727 \
v^8}{185515403673600}+\ldots
\nonumber \\
\ds b^5_{T,5}=\frac{50277247}{985600}-\frac{16301796103 \
v^2}{230630400}+\frac{263828905451 \
v^4}{5535129600}-\\
\ds \frac{400330948542829 \
v^6}{19760412672000}+\frac{54441481708734389 \
v^8}{9010748178432000}-\ldots
\nonumber \\
\ds b^5_{T,6}=-\frac{253491379}{4435200}+\frac{16301796103 \
v^2}{192192000}-\frac{117727186937 \
v^4}{1976832000}+\\
\ds \frac{3958504514434801 \
v^6}{148203095040000}-\frac{7164842887862063 \
v^8}{834328535040000}+\ldots
\nonumber
\ec
\clearpage
\section*{Figures}

\begin{figure}[hb]
\resizebox{1.00\textwidth}{!}
{\includegraphics{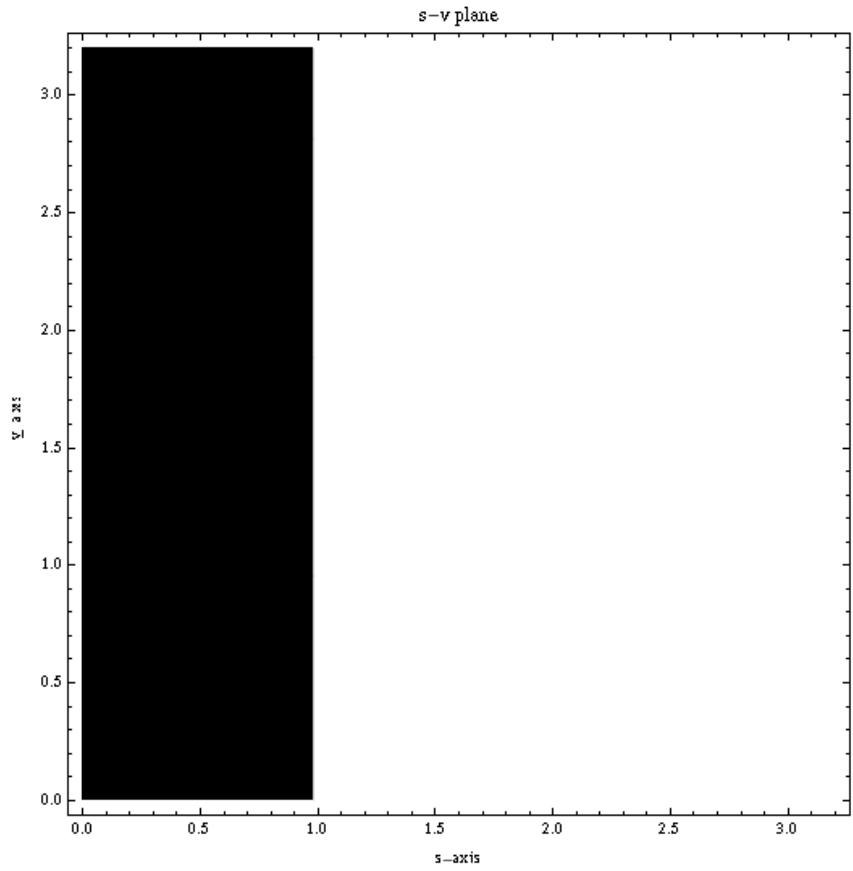}}
\caption{The stability region (s-v plane) of the classical Quinlan-Tremaine 12-step method}
\label{fig:1}
\end{figure}

\begin{figure*}
\resizebox{1.00\textwidth}{!}{%
\includegraphics{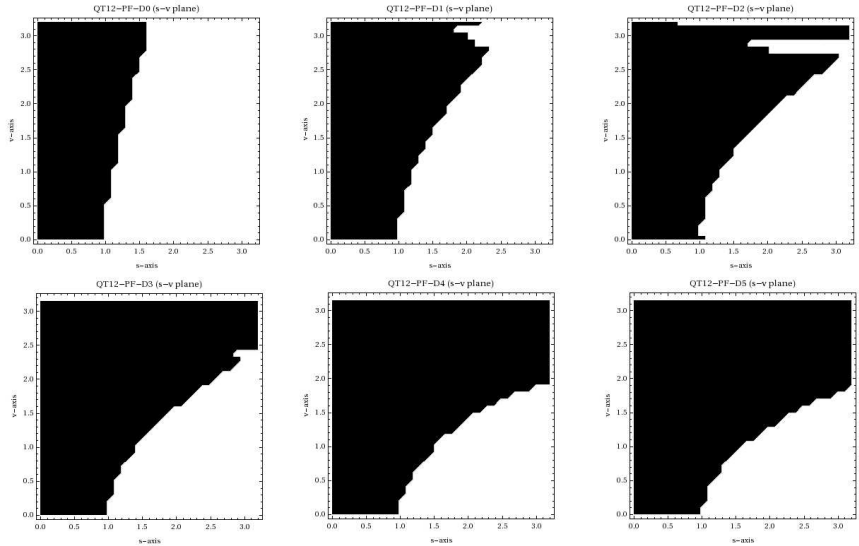}
}
\caption{The stability region of the methods PF-D0,PF-D1, PF-D2, PF-D3, PF-D4 and PF-D5 (from left to right and from top to bottom)}
\label{fig:2}
\end{figure*}

\begin{figure}
\resizebox{1.00\textwidth}{!}{%
\includegraphics{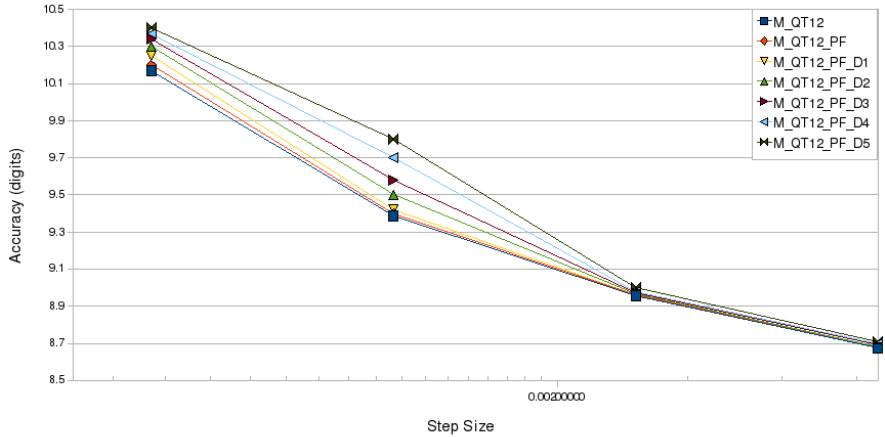}
}
\caption{The accuracy (digits) of the new methods compared to the classical one for the Schr\"odinger equation ($E=989.701916$)}
\label{fig:3}
\end{figure}

\begin{figure}
\resizebox{1.00\textwidth}{!}{%
\includegraphics{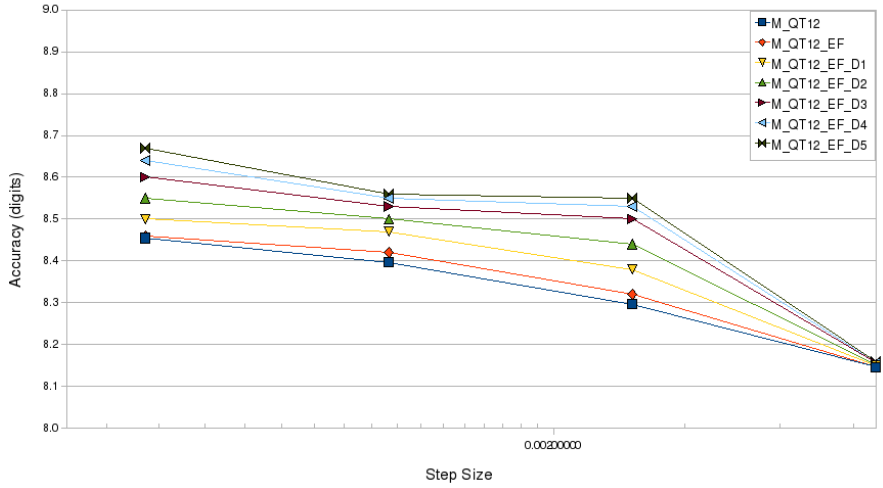}
}
\caption{The accuracy (digits) of the new methods compared to the classical one for the Schr\"odinger equation ($E=341.495874$)}
\label{fig:4}
\end{figure}

\begin{figure}
\resizebox{1.00\textwidth}{!}{%
\includegraphics{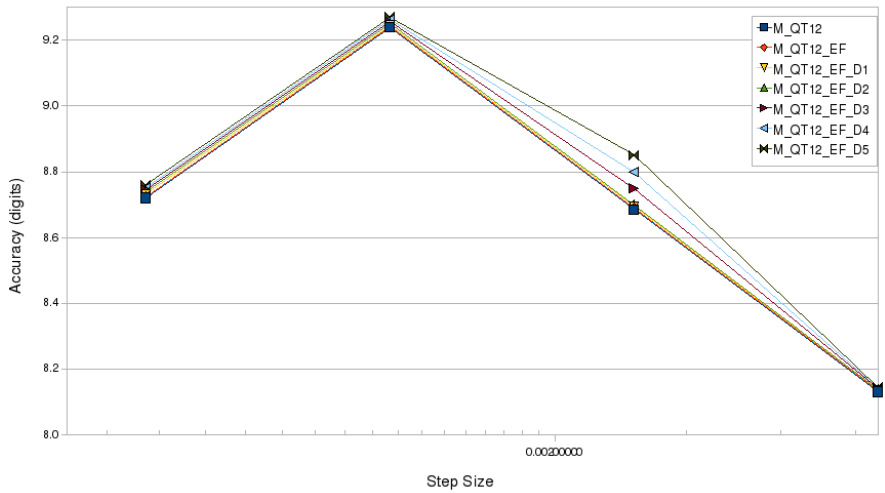}
}
\caption{The accuracy (digits) of the new methods compared to the classical one for the Schr\"odinger equation ($E=163.215341$)}
\label{fig:5}
\end{figure}

\end{article}
\end{document}